\documentclass[12pt]{article}
\usepackage{amssymb,amsmath,amsthm,abstract}
\usepackage{pifont,hyperref,geometry}
\usepackage{graphicx}
\usepackage{subcaption}
\usepackage{pdfpages}
\usepackage{pgfplots}
\usepackage{float}
\pgfplotsset{compat=1.18}
\usepackage{bm}
\usepackage{changepage}
\usepackage{setspace}
\usepackage{titlesec}
\usepackage{authblk}
\usepackage{fancyhdr}
\usepackage{xcolor}
\usepackage[numbers,sort]{natbib}
\titlespacing{\section}{0pt}{*0.3}{*0.3}
\newtheorem{lemma}{Lemma}[section]
\newtheorem{theorem}[lemma]{Theorem}
\newtheorem{definition}{Definition}[section]

\newtheorem{proposition}{Proposition}[section]
\newtheorem{remark}[lemma]{Remark}

\def\RR{{\mathbb{R}}}
\def\ZZ{{\mathbb{Z}}}
\def\TT{{\mathbb{T}}}

\def\D{{\mathcal{D}}}
\def\H{{\mathcal{H}}}

\def\d{{\mathrm{d}}}
\newcommand{\p}{\partial}
\newcommand{\nab}{\nabla}
\newcommand{\eps}{\varepsilon}

\title{\textbf{Small scale creation in 2D gravity-capillary water waves with vorticity}}
\author[1\thanks{E-mail: \texttt{typ23000101@mail.ustc.edu.cn}}]{Yuanpeng Tu}
\affil[1]{\small School of the Gifted Young, University of Science and Technology of China, Hefei, Anhui 230026, China.}
\date{\today}

\begin{document}
\numberwithin{equation}{section}
\maketitle

\begin{abstract} 
    In this paper, we consider 2D incompressible Euler equations in an unbounded domain with a free surface and a fixed bottom at finite depth. 
    The fluid motion is under the influence of gravity and surface tension. We construct initial data with a flat free surface and small velocity, such that the $L^\infty $ norm of the vorticity gradient has at least a double-exponential growth rate within the lifespan of the corresponding solution. 
    This work generalizes the result of Zlato\v{s} \cite{Zlatos2025} to the free-surface setting and Hu--Luo--Yao \cite{HLY2024} to the case of an unbounded domain.
\end{abstract}



\section{Introduction}

The 2D incompressible free boundary Euler equations describe the motion of an incompressible fluid in two dimensions with a free boundary separating the moving fluid region $\D_t $ and the vacuum region. In
the fluid region, the fluid velocity $u(t,x) $ and the pressure $p(t,x) $ satisfy the incompressible Euler
equations under the influence of gravity:
\begin{equation}\label{wave1}
    \begin{cases}
        \partial_t u+u\cdot\nabla u+\nabla p+ge_2=0, & \mbox{in } \D_t,\\
        \nabla\cdot u=0, & \mbox{in } \D_t,
    \end{cases}
\end{equation}
where $x=(x_1,x_2),\nab=(\p_1,\p_2)$ and $e_2=(0,1) $.
The fluid domain $\D_t $ lies between a fixed bottom $\Gamma_b=\RR\times\{x_2=0\} $ and an upper free surface $\Gamma_t $.
The free boundary $\Gamma_t $ evolves with the fluid velocity $u(t,x) $, namely, the speed of the free surface, denoted by $V$, is given by 
\begin{equation}\label{wave2}
    V=u\cdot N\quad \mbox{on }\Gamma_t,
\end{equation}
where $N $ is the outward unit normal to $\Gamma_t $. Throughout this paper, we assume the free-surface motion is under the influence of surface tension, and thus the pressure on the free boundary obeys
\begin{equation}\label{wave3}
    p=\sigma H\quad \mbox{on }\Gamma_t,
\end{equation}
where $\sigma > 0$ is the surface tension constant, and $H$ is the mean curvature of the free boundary.
On the fixed bottom $\Gamma_b$, we impose the slip boundary condition
\begin{equation}\label{wave4}
    u\cdot e_2=0\quad \mbox{on }\Gamma_b.
\end{equation}
For simplicity, we assume the initial free boundary to be a straight line $\Gamma_0 = \RR\times\{x_2 = 2\}$, so the initial fluid domain is
\begin{equation}\label{wave5}
    \D_0=\RR\times(0,2).
\end{equation}

\subsection{An overview of previous results}

The local well-posedness for the free-boundary incompressible Euler equations \eqref{wave1}--\eqref{wave4} has been well-studied in the past several decades. We refer to \cite{Wu1997,Wu1999,Lindblad2005,Lannes2005,ZZ2008,ABZ2014} for the case of no surface tension and \cite{CS2007,SZ1,SZ2,SZ3,MZ2009,ABZ2011} for the case of nonzero surface tension and references therein.
A natural further question is whether the small-data water-wave problem is globally well-posed. Positive answers to the long-time well-posedness have been given for the irrotational incompressible water waves with or without surface tension, and we refer to \cite{Wu2009,Wu2011,GMS2012,GMS2015,IP2015,IP2018,AD2015,IT2016,IT2017,DIPP2017,WangXC2018,WangXC2020,ZhengF2019} and references therein. For initial data with non-constant vorticity, there have been very few long-time existence results. Su \cite{Su2020} showed that, for 2D infinite-depth incompressible water waves with an initially downward-traveling symmetric vortex pair, the free interface remains smooth for an $O(\varepsilon^{-2})$-size lifespan. 

The goal of this paper is to give a negative answer to the question of long-time small-data stability of the water wave system with non-constant vorticity: for the 2D gravity-capillary water wave system with finite depth, we construct smooth initial data with a flat initial free surface and arbitrarily small velocity, such that the vorticity gradient $\|\nab\omega(t)\|_{L^\infty}$ admits a double-exponential growth rate within the lifespan of the solution. A straightforward corollary is that $\|u(t)\|_{W^{2,\infty}}$ grows to $O(1)$-size by the time scale $O(\varepsilon^{-1}\log\log \varepsilon^{-1})$ if the initial data that we construct is of size $\varepsilon\ll 1$, unless a singularity forms before this time.

As for the construction of fast growth for the vorticity gradient, a landmark breakthrough was established by Kiselev--\v{S}ver\`ak \cite{Sverak2013}, where they constructed smooth initial data for the 2D Euler equations in a disk such that the vorticity gradient has double-exponential growth, which was then extended to smooth domains with an axis of symmetry by Xu \cite{Xu2016}. The abovementioned works all require the fluid domain to be bounded, whereas a very recent breakthrough by Zlato\v{s} \cite{Zlatos2025} proves that double-exponential growth does occur for smooth solutions to the 2D Euler equations in the half-plane, an unbounded domain with a boundary. 

For the 2D Euler equations in a domain without boundary, Zlato\v{s} \cite{Zlatos2015} constructed a $C^{1,\alpha}$ solution in $\TT^2$ with exponential growth of the vorticity gradient. For smooth solutions, the first example of superlinear growth was constructed by Denisov \cite{Denisov2009}, and see also a recent result by Jeong--Yao--Zhou \cite{JYZ2025} about $O(t\log t)$-growth with the help of Lamb dipoles.

However, the small-scale creation in water waves is far less studied. To our knowledge, the only available result is due to Hu--Luo--Yao \cite{HLY2024}, which generalizes Kiselev--\v{S}ver\`ak \cite{Sverak2013} to the free-surface setting by constructing a smooth solution with double-exponential growth for the vorticity gradient for 2D periodic, finite-depth, incompressible water wave equations with surface tension, but the fluid domain considered in \cite{HLY2024} is still bounded. In this paper, we aim to tackle the case of an unbounded domain with a free surface, that is, the 2D incompressible gravity-capillary water wave system with finite depth. We aim to construct smooth initial data such that $\|\nabla\omega(t)\|_{L^\infty}$ grows at least double-exponentially for all times within the lifespan of the solution.

\subsection{Main results and strategies}
Our main result is as follows, which is stated for the $\sigma>0,g>0 $ case for simplicity.
\begin{theorem}\label{th1}
    For the 2D gravity-capillary water wave problem \eqref{wave1}--\eqref{wave4}, there exists a smooth velocity field $v_0\in C^\infty(\D_0)$ and universal constants $\varepsilon_0,c_1,c_2$, such that for any $\varepsilon\in(0,\varepsilon_0)$, 
    the solution\footnote{Here and throughout, a solution always means the velocity $u(t,\cdot)\in H^s(\D_t)$ and the domain $\D_t$ with $\Gamma_t\in H^{s+\frac32}$ that satisfy \eqref{wave1}--\eqref{wave4} for some fixed $s\geq 4 $.
    By elliptic regularity for the pressure associated with the free-boundary Euler system with boundary condition $p=\sigma H\in H^{s-\frac12} $ on $\Gamma_t $ and $\partial_{x_2} p=-g$ on $\Gamma_b $, we have $p(t,\cdot)\in H^s(\D_t) $. } 
    with initial velocity $u_0=\eps v_0$ and initial domain $\D_0$ given by \eqref{wave5} satisfies the following double-exponential growth rate for its vorticity gradient:
    \[\|\nabla\omega(t,\cdot)\|_{L^\infty(\D_t)}\geq\varepsilon\exp(c_1\exp(c_2\varepsilon t))\]
    for all $t\in[0,T) $, where $T $ denotes the lifespan of the solution. Here $\omega(t,x):=\p_1 u_2-\p_2 u_1$ is the vorticity of the velocity $u(t,x)$.
\end{theorem}
\begin{remark}
    A straightforward corollary of this theorem is that we construct $C^\infty$ small data of size $\eps(\ll 1)$ such that $\|u\|_{W^{2,\infty}}$ grows to $O(1)$-size by the time scale $O(\eps^{-1}\log\log\eps^{-1})$ unless some singularity forms before this time. 
    This corollary is a sharp contrast to the long-time existence results for the irrotational case. 
\end{remark}

To prove Theorem \ref{th1}, a natural starting point is to impose the same symmetry as \cite{HLY2024}, with $u_{01} $ odd-in-$x_1 $ and $u_{02} $ even-in-$x_1 $ respectively.
In fact, this symmetry holds throughout the lifespan, so the vorticity remains odd-in-$x_1 $ for all time.
If we set $\omega_0=\varepsilon $ at some special points in the right half of $\D_0 $, as in \cite{Zlatos2025}, it can be shown that there exists a pair of small regions that propel themselves toward each other at a double-exponential rate, which completely fill a shrinking neighborhood of the origin.
These regions have distinct vorticity values, touch the boundary of the half plane and thus lead to the double-exponential growth for the vorticity gradient near the origin on the fixed bottom. 
See Figure \ref{F1} for a visualized explanation.
\begin{figure}[H]
    \centering 
    \includegraphics[scale=0.45]{snapshots.png}
    \caption{Illustrations of our initial data $\omega_0,\D_0 $, together with their evolution at later times. The four panels correspond to $t=0,1,2,3$ (from left to right and top to bottom). Here the blue and red colors represent positive and negative vorticity respectively; the black and red curves denote the free surface and the fixed bottom. As we show in Proposition \ref{pp2}, the free boundary $\Gamma_t $ remains close to $\Gamma_0$ throughout the lifespan of the solution.} 
    \label{F1}
\end{figure}

The main challenge lies in the deformation of the domain.
Because $\omega\neq 0 $ and $\D_t $ evolves in time, the fluid domain $\D_t $ may deviate significantly from the initial domain $\D_0 $.
In particular, the free boundary $\Gamma_t $ might approach arbitrarily close to the origin within finite time. Motivated by \cite{HLY2024}, we show in Proposition \ref{pp2} that for small initial data, this never occurs within the lifespan of the solution by using the energy conservation. Indeed, the free boundary Euler equations with surface tension possess a conserved energy $E(t)=K(t)+P(t)$, where $K(t)$ denotes the kinetic energy and $P(t)$ denotes the potential energy. 
More precisely, the surface-tension part in $P(t) $ controls the excess length of the free boundary, while the gravitational part in $P(t) $ controls the vertical displacement of the fluid region.
This would force the free surface lie in a fixed neighborhood of the initial flat surface $\{x_2=2\} $ away from the bottom, provided $K(0) $ is sufficiently small.

It is straightforward to prove this when $\Gamma_t$ is the graph of a function.
For the general case, Hu--Luo--Yao \cite{HLY2024} observed that the potential energy is proportional to the perimeter of the free boundary $\Gamma_t$ for 2D periodic water waves with surface tension. However, in our case, the horizontal component is $\RR$ instead of $\TT$, and therefore the perimeter of $\Gamma_t$ is always infinite.

To overcome this difficulty arising from the unboundedness of the domain, we shall appropriately employ Steiner symmetrization to convert a ``non-graph'' free boundary into a ``graph'' boundary. Since the Steiner symmetrization is only useful for bounded domains, we restrict attention to the portion of $\D_t $ contained in $\{|x_1|<R\}$. 
We apply Steiner symmetrization twice (with respect to two different directions) to both the domain above the initial free boundary $\D_+(R):=\D_t(R)\setminus \D_0$ and the domain under the initial fluid $\D_-(R):=(\D_0\setminus\D_t(R))\cap\{|x_1|<R\}$. 
The first Steiner symmetrization $S_{e_1}$ relocates the highest and lowest points within the truncated domains to $x_1=0$, and then the second Steiner symmetrization $S_{e_2}$ transforms these irregular regions into ``regular'' regions whose boundaries are graphs of some functions.
By relocating the highest and lowest points, we guarantee that their relative elevations are preserved under the second Steiner symmetrization.

Another difficulty is that the exact Biot-Savart law is not available in the free boundary setting. However, the approximate Biot-Savart law in \cite{HLY2024} allows us to obtain a pointwise estimate of $u $ analogous to the key lemma \cite[Lemma 3.1]{Sverak2013}, leading to the double-exponential growth of $\|\nabla\omega(t,\cdot)\|_{L^\infty(\D_t)} $. 

\paragraph*{Acknowledgment.} The author would like to thank Prof. Junyan Zhang, the mentor of the author's undergraduate research project, for suggesting this problem and revising the draft version of this manuscript.

\section{Preliminary results}
This section records some preliminary results that are important in the proof of Theorem \ref{th1}, 
including the symmetry of 2D incompressible free-surface Euler equations, the energy conservation, the separation of the free surface, and the fixed bottom and key estimates for the fluid velocity and the error of the approximate Biot--Savart law.

\subsection{Symmetry of 2D free-surface Euler equations}
To begin with, we discuss some symmetry properties of the 2D free boundary Euler equations. 
We state the symmetry assumption in terms of the velocity rather than the vorticity since one cannot uniquely determine the velocity using the vorticity at a given moment.
\begin{lemma}\label{sym}
    Let $(u_0,\D_0) $ be the initial data of \eqref{wave1}--\eqref{wave4} with initial domain $\D_0 $ given by $\eqref{wave5}$ and initial velocity $u_0=(u_{01},u_{02}) $ satisfying
    \[-u_{01}(-x_1,x_2)=u_{01}(x_1,x_2),\quad u_{02}(-x_1,x_2)=u_{02}(x_1,x_2).\]
    Then for all times during the lifespan of the solution, the solution $(u,\D_t) $ satisfies the same symmetry, that is
    \[-u_1(t,-x_1,x_2)=u_1(t,x_1,x_2),\quad u_2(t,-x_1,x_2)=u_2(t,x_1,x_2),\]
    and the moving fluid domain $\D_t $ remains even in $x_1 $, that is 
    \[\D_t=\tilde{\D}_t=\{(-x_1,x_2):(x_1,x_2)\in\D_t\}.\]
\end{lemma}
\begin{proof}
See Hu--Luo--Yao \cite[Lemma 2.1]{HLY2024}.
\end{proof}
\begin{remark}
    As a direct consequence of Lemma \ref{sym}, we have that the vorticity $\omega(t,x)=\nabla^\bot\cdot u(t,x) $ stays odd in $x_1 $ during the lifespan of the solution.
\end{remark}

\subsection{Conservation of energy}
In this subsection, we aim to prove that the free surface $\Gamma_t$ never touches the fixed bottom $\Gamma_b$ within the lifespan of the solution. To achieve this, we shall first prove the energy conservation.
We consider a renormalized energy functional to compensate for the infinite horizontal extent.
\begin{proposition}\label{prop1}
    Take $x_2=1 $ as the zero level of the gravitational potential energy.
    Let $(u,\D_t) $ be a solution to the system \eqref{wave1}--\eqref{wave4}. Define the renormalized energy
    \begin{equation}
        E(t)=K(t)+P(t),
    \end{equation}
    where 
    \[K(t):=\frac12\int_{\D_t}{|u(t,x)|}^2\,dx\]
    and
    \[P(t):= \sigma \lim_{R \to \infty} \left( \int_{\Gamma_t \cap \{|x_1| < R\}} \,d S_t - 2R \right) +\lim_{R\to\infty}\int_{\D_t\cap\{|x_1|<R\}} g(x_2-1)\,dx\]
    are the kinetic energy and the potential energy respectively. Then the above limits are finite, and the energy is conserved:
    \begin{equation}\label{ene}
        E(t)=E(0)
    \end{equation}
    for all times within the lifespan of the solution.
\end{proposition}
\begin{proof}
For $R>0$, define
\[\D_t(R):=\D_t\cap\{-R<x_1<R\},\qquad\Gamma_t(R):=\Gamma_t\cap\{-R<x_1<R\},\]
and let
\[\Sigma_t^{\pm}(R):=\D_t\cap\{x_1=\pm R\}\]
denote the two vertical lateral sides of $\partial \D_t(R)$. Set $E_R(t,R)=K_R(t)+P_R(t) $, where
\[K_R(t):=\frac12\int_{\D_t(R)}|u(t,x)|^2\,dx\]
and
\[P_R(t):=\sigma\left(\int_{\Gamma_t(R)}dS_t-2R\right)+\int_{\D_t(R)}g(x_2-1)\,dx.\]

\medskip
\noindent
\textbf{Step 1. Time derivative of $K_R(t)$.}
By Reynolds' transport formula and the fact that the free surface moves with normal velocity $u\cdot N$, while the bottom and the two vertical
sides are fixed, we have
\[\frac{d}{dt}K_R(t)=\int_{\D_t(R)}u\cdot \partial_t u\,dx+\frac12\int_{\Gamma_t(R)}|u|^2(u\cdot N)\,dS_t.\]
On the other hand, using $\nabla\cdot u=0$ and the divergence theorem, we obtain
\[\int_{\D_t(R)}(u\cdot\nabla u)\cdot u\,dx=\int_{\partial \D_t(R)}\frac12|u|^2(u\cdot N)\,dS_t.\]
Since $u\cdot e_2=0$ on the fixed bottom $\Gamma_b$, this becomes
\[\int_{\D_t(R)}(u\cdot\nabla u)\cdot u\,dx=\frac12\int_{\Gamma_t(R)}|u|^2(u\cdot N)\,dS_t+\frac12\int_{\Sigma_t^+(R)\cup\Sigma_t^-(R)}|u|^2(u\cdot N)\,dS_t.\]
Therefore,
\begin{align*}
    \frac{d}{dt}K_R(t)&=\int_{\D_t(R)}u\cdot(\partial_t u+u\cdot\nabla u)\,dx-\frac12\int_{\Sigma_t^+(R)\cup\Sigma_t^-(R)}|u|^2(u\cdot N)\,dS_t.\\
    &=-\int_{\D_t(R)}\nabla p\cdot u\,dx-\int_{\D_t(R)}g\,u_2\,dx-\frac12\int_{\Sigma_t^+(R)\cup\Sigma_t^-(R)}|u|^2(u\cdot N)\,dS_t.
\end{align*}
Applying the divergence theorem again and using $\nabla\cdot u=0$, we get
\[-\int_{\D_t(R)}\nabla p\cdot u\,dx=-\int_{\partial \D_t(R)}p(u\cdot N)\,dS_t=
-\int_{\Gamma_t(R)}p(u\cdot N)\,dS_t-\int_{\Sigma_t^+(R)\cup\Sigma_t^-(R)}p(u\cdot N)\,dS_t,\]
Similarly,
\begin{align*}
    -\int_{\D_t(R)}g\,u_2\,dx&= -\int_{\D_t(R)}\nabla\cdot\bigl(g(x_2-1)u\bigr)\,dx\\
    &=-\int_{\Gamma_t(R)}g(x_2-1)(u\cdot N)\,dS_t-\int_{\Sigma_t^+(R)\cup\Sigma_t^-(R)}g(x_2-1)(u\cdot N)\,dS_t.
\end{align*}
Since $p=\sigma H$ and $V=u\cdot N$ on $\Gamma_t$,
\begin{align*}
    \frac{d}{dt}K_R(t)&=-\int_{\Gamma_t(R)}\bigl(p+g(x_2-1)\bigr)(u\cdot N)\,dS_t+\mathcal R_K(R,t)\\
    &=-\int_{\Gamma_t(R)}\sigma H V\,dS_t-\int_{\Gamma_t(R)}g(x_2-1)V\,dS_t+\mathcal R_K(R,t),
\end{align*}
where
\begin{equation}\label{2.3}
\mathcal R_K(R,t):=-\int_{\Sigma_t^+(R)\cup\Sigma_t^-(R)}\left(p(u\cdot N)+g(x_2-1)(u\cdot N)+\frac12|u|^2(u\cdot N)\right)\,dS_t.
\end{equation}

\medskip
\noindent
\textbf{Step 2. Time derivative of $P_R(t)$.}
By the transport formula for an evolving hypersurface, the area element satisfies
\[\frac{d}{dt}\int_{\Gamma_t(R)}dS_t=\int_{\Gamma_t(R)}H V\,dS_t+\mathcal R_\sigma(R,t),\]
where $\mathcal R_\sigma(R,t)={\left[\tau\cdot u\right]}_{\partial \Gamma_t(R)}$ with $\tau $ is the unit tangent vector of $\Gamma_t $ which represents the boundary contribution from the endpoints of $\Gamma_t $.
Hence
\begin{equation}\label{2.4}
    \frac{d}{dt}\left[\sigma\left(\int_{\Gamma_t(R)}dS_t-2R\right)\right]=\int_{\Gamma_t(R)}\sigma H V\,dS_t+\sigma\mathcal R_\sigma(R,t). 
\end{equation}
For the gravitational potential energy, Reynolds' transport formula gives
\begin{equation}\label{2.5}
    \frac{d}{dt}\int_{\D_t(R)}g(x_2-1)\,dx=\int_{\Gamma_t(R)}g(x_2-1)V\,dS_t 
\end{equation}
since the bottom and the two vertical sides are fixed.

\medskip
\noindent
\textbf{Step 3. Cancellation and limit.}
Combining \eqref{2.3}, \eqref{2.4}, and \eqref{2.5}, we obtain
\[\frac{d}{dt}E_R(t,R)=\mathcal R_K(R,t)+\sigma\mathcal R_\sigma(R,t) .\]
Since $u(t,\cdot),p(t,\cdot)\in H^s(\D_t)$ with $s\ge 4$,
we have $u\in C^1$ and $u|_{\Gamma_t}\in H^{s-\frac12}(\Gamma_t) $.
Furthermore, it holds that $\int_{\Sigma_t^\pm(R)} |u|^2\,dS_t $ is $W^{1,1}(\RR) $ function, taking $R $ as variable. 
Hence, for every compact interval $[0,T_0]\subset [0,T) $, there exists a constant $M>2 $
such that for any $R $, $\Sigma_t^\pm(R)\subseteq \{\pm R\}\times (0,M) $, and as $R \to\infty $
\[\int_{\Sigma_t^\pm(R)} |u|^2\,dS_t \to 0\quad\mbox{and}\quad \|u\|_{L^\infty(\partial \Gamma_t(R))}\to 0\]
for any $t\in [0,T_0] $.
Sobolev embedding implies that $u$ and $p$ are bounded on each time slice and continuous in time. We obtain
\begin{align*}
    |\mathcal R_K(R,t)|&\leq \int_{\Sigma_t^+(R)\cup\Sigma_t^-(R)}\left(p|u|+g |x_2-1||u|+\frac12|u|^3\right)\,dS_t\\
    & \leq \left(\|p(t)\|_{L^\infty}+gM\right)\sqrt{M\int_{\Sigma_t^+(R)\cup\Sigma_t^-(R)}{|u|}^2\,dS_t}+\frac12\|u(t)\|_{L^\infty}\int_{\Sigma_t^+(R)\cup\Sigma_t^-(R)}{|u|}^2\,dS_t.
\end{align*}
Therefore, as $R\to\infty $,
\[\mathcal R_K(R,t)+\sigma\mathcal R_\sigma(R,t) \to 0\]
uniformly in $t\in [0,T_0] $.
Since $E(0)=\frac12\int_{\D_0}{|u_0|}^2\,dx<\infty $, 
\[E_R(t,R)=E_R(0,R)+\int^t_0\mathcal R_K(R,s)+\sigma\mathcal R_\sigma(R,s)\d s. \]
Passing to the limit $R\to\infty$ yields for any $t\in [0,T_0] $,
\[E(t)=\lim_{R\to\infty}E_R(t,R)=\lim_{R\to\infty}\left(E_R(0,R)+\int^t_0\mathcal R_K(R,s)+\sigma\mathcal R_\sigma(R,s)\d s\right)=E(0). \]
Since $T_0<T$ is arbitrary, \eqref{ene} holds throughout the lifespan of the solution.
\end{proof}

Using energy conservation, we prove that — provided the initial free boundary is flat and the initial kinetic energy is sufficiently small — the free surface never touches the bottom, a fact that lies at the heart of the proof of the paper's main theorem.
\begin{proposition}\label{pp2}
    Under the assumption of Theorem \ref{th1}, we additionally assume that there exist two constants $C_1>0$ and $\varepsilon>0 $, with $\varepsilon $ sufficiently small, such that $K(0)=E(0)<C_1\varepsilon^2 $. Then we have 
    \begin{equation}\label{p2}
        \Gamma_t\subset \RR\times\left(\frac32,\frac52\right)\quad \forall\, t\in [0,T).
    \end{equation}
    In particular, this shows that $\Lambda =\RR\times[0,1]\subset \D_t $ within the lifespan of the solution. Moreover,
    \begin{equation}\label{p3}
        K(t)\leq K(0) \quad \forall\, t\in [0,T).
    \end{equation}
\end{proposition}
\begin{proof}
    Assume, to the contrary, that there exists $t\in(0,T) $ such that $\Gamma_t\not\subset \RR\times\left(\frac32,\frac52\right) $.
    
    \textbf{Case 1: The free surface is a graph within the lifespan.} 
    We also assume that free surface is closed to $x_2=2 $ when $|x_1| $ is large.
    In such a case, we can prove the result by straightforward energy estimates.
    
    Assume $\Gamma_t=\{(x_1,2+\varphi(x_1)):x_1\in\RR\} $. By conservation of volume, we have 
    \[\|\varphi\|_{L^\infty(\RR)}\geq \frac12\qquad\mbox{and}\qquad \int_\RR\varphi(x_1)\,dx_1=0, \]
    and $P(t)=P_1(\varphi)+P_2(\varphi) $, where 
    \[P_1(\varphi)=\int_{\RR}\sigma\left(\sqrt{1+{\varphi'(x_1)}^2}-1\right)\,dx_1\]
    and 
    \[P_2(\varphi)=\int_{\RR}g\frac{{(\varphi(x_1)+1)}^2-1}2\,dx_1=\frac{g}2\int_{\RR}{\varphi(x_1)}^2\,dx_1.\]
    Notice that 
    \[\sqrt{1+{\varphi'(x_1)}^2}-1\geq \frac1{1+\sqrt{2}}\min\left\{|\varphi'(x_1)|,{(\varphi'(x_1))}^2\right\}.\]
    Without loss of generality, assume $\varphi(0)=\frac12,\varphi(R)=\frac14 $, and $\varphi(x_1)\in [\frac14,\frac12] $ for all $x_1\in [0,R] $.
    Let $A=\{x_1:|\varphi'(x_1)|\geq 1\}\cap [0,R] $ and $B=\{x_1:|\varphi'(x_1)|< 1\}\cap [0,R] $.
    Then we have: 
    \begin{align*}
        \frac14& \leq \int_{[0,R]}|\varphi'(x_1)|\,dx_1=\int_A|\varphi'(x_1)|\,dx_1+\int_B|\varphi'(x_1)|\,dx_1\\
        & \leq \int^R_0\min\left\{|\varphi'(x_1)|,{(\varphi'(x_1))}^2\right\}\,dx_1+\sqrt{R\int^R_0\min\left\{|\varphi'(x_1)|,{(\varphi'(x_1))}^2\right\}\,dx_1}\\
        & \leq \frac32\int^R_0\min\left\{|\varphi'(x_1)|,{(\varphi'(x_1))}^2\right\}\,dx_1+\frac{R}{2}\\
        & \leq \frac{3(1+\sqrt{2})}{2\sigma}P_1(\varphi)+\frac{R}{2} .
    \end{align*}
    Since $\varphi^2(x_1)\geq 1/16 $ in $[0,R] $, it follows that  
    \[P_2(\varphi)=\frac{g}2\int_{\RR}{\varphi(x_1)}^2\,dx_1\geq \frac{Rg}{32}.\]
    Consequently, we obtain 
    \[P(t)\cdot {4\min\left\{\frac{3(1+\sqrt{2})}{2\sigma},\frac{16}{g}\right\}}\geq 1. \]
    This yields $P(t)\geq c>0 $, contradicting $E(t)=E(0)\leq C_1\varepsilon^2 $ for $\varepsilon $ small.
    Thus, if we take $\varepsilon $ sufficiently small, then $\eqref{p2}$ holds.
    Then $\eqref{p3}$ follows from the energy conservation and the fact that $P_1(\varphi),P_2(\varphi)\geq 0 $.

    \textbf{Case 2: General case.} Now $\Gamma_t$ is not necessarily a graph over $x_1 $, but we can still  convert this case to the graph case with the help of Steiner symmetrization.
    \begin{definition}[Steiner symmetrization]\label{S}
        Given $a,b\in\RR^2 $, with $|a|=1 $. We denote $L^a_b=\{b+ta:t\in\RR\} $ to be the line through $b $ in the direction $a $, and $P_a=\{x\in\RR^2:x\cdot a=0\} $ to be the plane through the origin perpendicular to $a $.
        We define the Steiner symmetrization of domain $A $ with respect to the plane $P_a $ to be the set 
        \[S_a(A)=\bigcup_{\substack{b\in P_a\\ A\cap L^a_b\neq\emptyset}}\left\{b+ta:|t|\leq\frac12\H^1(A\cap L_b^a)\right\}\]
        where $\H^1 $ is the 1D Hausdorff measure.
    \end{definition}
    For any $\delta>0 $, choose sufficiently large constant $R$ such that 
    \[\{\pm R\}\times(0,2-\delta]\subseteq\Sigma_t^{\pm}(R)\subseteq \{\pm R\}\times(0,2+\delta)\]
    and  
    \[\left|\int_{\Gamma_t\setminus\Gamma_t(R)}\sigma\,d (S_t-x_1)+\lim_{M\to\infty}\int_{\D_t(M)\setminus\D_t(R)}g(x_2-1)\,dx\right|<\frac1{5}E(0),\]
    where $\D_t(R)=\D_t\cap\{-R<x_1<R\} $, $\Gamma_t(R)=\Gamma_t\cap\{-R<x_1<R\} $ and $\Sigma_t^{\pm}(R):=\D_t\cap\{x_1=\pm R\} $.

    Consider the domain above the initial fluid $\D_+(R)=\D_t(R)\setminus \D_0$ and the domain under the initial fluid $\D_-(R)=(\D_0\setminus\D_t(R))\cap \{|x_1|<R\} $.
    Although the zero gravitational potential level is $x_2=1 $ , it is convenient to rewrite the gravitational contribution relative to the reference interface $x_2=2 $ after decomposing the domain into the regions above and below $\D_0 $. 
    We can rewrite the potential energy in the strip $-R<x_1<R $ as
    \begin{align*}
        P_R(t) =& \int_{\Gamma_t(R)} \sigma\,d (S_t-x_1)  +\int_{\D_t(R)}g(x_2-1)\,dx\\
        = & \int_{\Gamma_+(R)\cup \Gamma_-(R)} \sigma\,d (S_t-x_1)+\int_{\D_+(R)}g(x_2-2)\,dx + \int_{\D_-(R)}g(2-x_2)\,dx.
    \end{align*}
    where $\Gamma_+(R)=\Gamma_t(R)\cap\overline{\D_+(R)} $ is the free boundary of $\D_+(R) $ and $\Gamma_-(R)=\Gamma_t(R)\cap\overline{\D_-(R)} $ is the free boundary of $\D_-(R) $.

    We define two domains that are symmetric with respect to the $x_1$-axis:
    \[\tilde{\D}_+(R):=\{(x_1,x_2)\in \RR^2:(x_1,2+|x_2|)\in \D_+(R)\} \]
    and 
    \[\tilde{\D}_-(R):=\{(x_1,x_2)\in \RR^2:(x_1,2-|x_2|)\in \D_-(R)\}.\]
    Consider the transformations of these domains
    \[\D^*_+(R):=S_{e_2}(S_{e_1}(\tilde{\D}_+(R)))\cap (\RR\times\RR_+),\]
    \[\D^*_-(R):=S_{e_2}(S_{e_1}(\tilde{\D}_-(R)))\cap (\RR\times\RR_+).\]
    where $S $ denotes the Steiner symmetrization introduced in Definition \ref{S}.
    See Figure \ref{F2} for an illustration.
    \begin{figure}[htbp]
        \centering
        \includegraphics[scale=0.65]{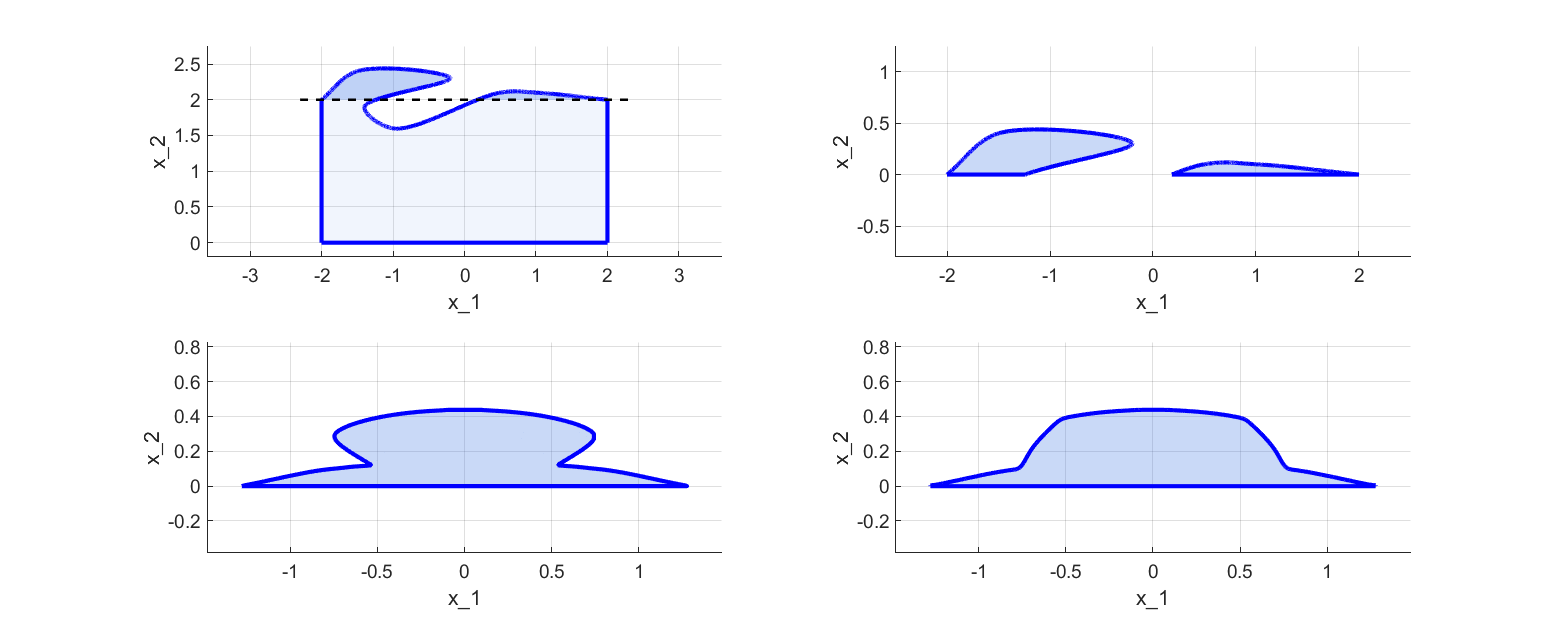}
        \caption{The first figure illustrates an example of the truncated domain $\D_t(R) $, where $R=2 $ in this example. The second figure depicts $\tilde{\D}_+(R) $ above $x_2=0 $. The third figure depicts the domain $S_{e_1}(\tilde{\D}_+(R)) $ above $x_2=0 $. The fourth figure depicts the domain $S_{e_2}(S_{e_1}(\tilde{\D}_+(R)))$ above $x_2=0 $.}
        \label{F2}
    \end{figure}

    After applying Steiner symmetrization, 
    the perimeter of the new graph does not increase, due to the isoperimetric inequality of Ennio De Giorgi \cite{Steiner},
    which means that 
    \[\int_{\Gamma_+^*(R)}\sigma\,d (S_t-x_1)\leq \int_{\Gamma_+(R)}\sigma\,d (S_t-x_1)+2\delta\sigma,\]
    and 
    \[\int_{\Gamma_-^*(R)}\sigma\,d (S_t-x_1)\leq \int_{\Gamma_-(R)}\sigma\,d (S_t-x_1)+2\delta\sigma,\]
    where $\Gamma^*_+(R) $ is the portion of the boundary of $\D^*_+(R) $ that lies above the $x_1$-axis and $\Gamma^*_-(R) $ is the portion of the boundary of $\D^*_-(R) $ that lies above the $x_1$-axis.
    Here the error term $2\delta\sigma$ comes from the two lateral truncation segments at
    $x_1=\pm R$: by the choice of $R$, Since $\{\pm R\}\times(0,2-\delta]\subseteq\Sigma_t^{\pm}(R)\subseteq \{\pm R\}\times(0,2+\delta) $, 
    So \[\partial \D_+(R)\cap\{x_1=\pm R\},\partial \D_-(R)\cap\{x_1=\pm R\}\subseteq \{\pm R\}\times [0,\delta].\]
    It shows that each truncated side contributes at most $\delta$ to the excess surface length above the flat reference configuration.
    
    Recall that in the preceding Steiner symmetrization procedure, we take the interface \(x_2=2\) as the reference level and translate the parts of \(\D_+(R)\) and \(\D_-(R)\) above and below this interface so that, after translation, the reference level \(x_2=2\) is moved to \(x_2=0\). 
    Thus, when comparing gravitational potential energies for the transformed domains \(\D^*_+(R)\) and \(\D^*_-(R)\), it is natural to regard \(x_2=0\) as the zero potential level for comparing the gravitational terms.
    Under this convention, $S_{e_1} $ does not change the gravitational potential energy in $\D_+(R) $ above the $x_1$-axis, since it only rearranges the set horizontally without affecting the \(x_2\)-coordinate.
    On the other hand, $S_{e_2} $ decreases the gravitational potential energy in $\D_+(R) $ above the $x_1$-axis.
    Indeed, this symmetrization rearranges each vertical slice into a segment symmetric with respect to the line \(x_2=0\), thereby moving mass from higher values of \(x_2\) to lower ones and reducing the average height.
    Therefore, the gravitational potential energy in $\D_+(R)$ is larger than that in $\D^*_+(R) $:
    \[\int_{\D^*_+(R)}g x_2\,dx\leq \int_{\D_+(R)}g(x_2-2)\,dx.\]
    In the same way, we obtain
    \[\int_{\D^*_-(R)}g x_2\,dx\leq \int_{\D_-(R)}g(2-x_2)\,dx.\]
    We thus get the estimate 
    \begin{align}\label{P}
        P_R(t) =& \sigma \left( \int_{\Gamma_+(R)\cup \Gamma_-(R)} \,d (S_t-x_1)\right) +\int_{\D_+(R)}g(x_2-2)\,dx + \int_{\D_-(R)}g(2-x_2)\,dx \nonumber\\
        \geq &\int_{\Gamma_+^*(R)}\sigma\,d (S_t-x_1) + \int_{\Gamma_-^*(R)}\sigma\,d (S_t-x_1) + \int_{\D^*_+(R)}g x_2\,dx + \int_{\D^*_-(R)}g x_2\,dx-4\delta\sigma . 
    \end{align}

    Since $\Gamma_t\not\subset \RR\times\left(3/2,5/2\right) $, it follows that either 
    \[\Gamma_+(R)\not\subset \RR\times\left(2,5/2\right) \quad\mbox{or}\quad \Gamma_-(R)\not\subset \RR\times\left(3/2,2\right) .\]
    Without loss of generality, assume that $\Gamma_-(R)\not\subset \RR\times\left(3/2,2\right) $. 
    From the construction of $\Gamma_-^*(R) $ and the following observation, it follows that $\Gamma^*_-(R)\not\subset \RR\times\left(0,1/2\right) $. Let 
    \[x^*_2:=\sup\{x_2>0:(x_1,x_2)\in \tilde{\D}_-(R)\}> 1/2 .\]
    Now apply the first Steiner symmetrization $S_{e_1}$. Since $S_{e_1}$ rearranges each horizontal slice only in the $x_1$-direction, it preserves, for every fixed $s$, whether the slice at height $x_2=s$ is empty or not. Hence
    \[\sup\{x_2>0:(x_1,x_2)\in S_{e_1}(\tilde{\D}_-(R))\}=x^*_2> 1/2 .\]
    Moreover, because $\tilde{\D}_-(R)$ is symmetric with respect to the line $x_2=0$, the set
    $S_{e_1}(\tilde{\D}_-(R))$ is still symmetric with respect to $x_2=0$.
    By the definition of \(x_2^*\), these slices are nonempty for all \(0<s<x_2^*\).
    Therefore, 
    \[\{0\}\times (-x^*_2,x^*_2) \subset S_{e_1}(\tilde{\D}_-(R)).\]
    Next we apply the second Steiner symmetrization $S_{e_2}$. This symmetrization acts on
    each vertical line $\{x_1=\mathrm{const}\}$ by replacing the intersection with a centered
    vertical interval of the same length.
    In particular, the vertical section of $S_{e_2}(S_{e_1}(\tilde{\D}_-(R)))$ at $x_1=0$ is exactly a centered interval containing
    $\{0\}\times(-x_2^*,x_2^*)$. After restricting to the upper half-plane, we conclude that 
    \[\Gamma_-^*(R)\not\subset \RR\times(0,1/2).\] 

    Finally, by construction of the second Steiner symmetrization, for each fixed $x_1$ the
    vertical section of $\D_-^*(R)$ is either empty or an interval of the form $[0,\varphi(x_1))$.
    Therefore $\Gamma_-^*(R)$ is the graph of a function $x_2=\varphi(x_1)$. By combining the proof of Case 1 with inequality \eqref{P}, we may take $\varepsilon,\delta $ sufficiently small so that 
    \[P_R(t)\geq \int_{\Gamma_-^*(R)}\sigma\,d (S_t-x_1)+\int_{\D^*_-(R)}g x_2\,dx-4\delta\sigma > \frac{6}{5}E(0).\]
    The following estimate then holds:
    \begin{align*}
        P(t)= & \lim_{R'\to\infty}\int_{\Gamma_t(R')\setminus\Gamma_t(R)}\sigma\,d (S_t-x_1)+\lim_{R'\to\infty}\int_{\D_t(R')\setminus\D_t(R)}g(x_2-1)\,dx+ P_R(t)\\
        > &-\frac15 E(0)+\frac{6}{5}E(0)=E(0) .
    \end{align*}
    The inequality $E(t)\geq P(t)>E(0) $ contradicts the conservation of energy stated in Proposition \ref{prop1}. 
    Hence, \eqref{p2} must hold throughout the lifespan.

    By the same argument, one can also prove $P(t)\geq 0 $. 
    Consequently, \eqref{p3} follows from the conservation of energy.
\end{proof}
\begin{remark}
    In fact, from the argument of Proposition \ref{pp2}, we can get the decay of $\Gamma_t $ when $|x_1| $ is large, that is for any $\delta>0 $, there exists $R>0 $ such that 
    \[\Gamma_t\setminus\Gamma_t(R)\subseteq \RR\times (2-\delta,2+\delta).\]
\end{remark}

\subsection{Error estimates of an approximate Biot-Savart law}
Because of the lack of the Biot-Savart law in the free boundary setting, we introduce an approximate Biot–Savart law from \cite[Proposition 3.2]{HLY2024} which only uses the information of $\omega(t,\cdot) $ in the set $\Lambda =\RR\times[0,1] $.
It provides a precise estimate of the error between the actual velocity $u $ and the approximate velocity $U $ near the origin.

For any $t\geq 0$ during the lifespan of the solution, we define an approximate velocity field $U(t,\cdot):\Lambda \to\RR^2 $
\begin{equation}
    U(t,\cdot)=\nabla^{\bot}\Psi(t,\cdot)\quad\mbox{in }\Lambda , 
\end{equation}
where $\Psi(t,\cdot) $ solves the following elliptic equation at the fixed time $t $:
\begin{equation}\label{Psi}
    \begin{cases}
        \Delta\Psi(t,\cdot)=\omega(t,\cdot) & \mbox{in }\Lambda \\
        \Psi(t,\cdot)=0 & \mbox{on }\partial\Lambda ,
    \end{cases}
\end{equation}
where $\omega(t,\cdot) $ is the vorticity of the solution $u(t,\cdot) $. We also define the error $e(t,\cdot):\Lambda \to\RR^2 $ as 
\begin{equation}\label{error}
    e(t,\cdot)=u(t,\cdot)|_\Lambda -U(t,\cdot).
\end{equation}

\begin{lemma}\label{err}
    Under the assumption of Proposition $\ref{pp2}$, $e(t,\cdot) $ is smooth in $B_{1/2}\cap\Lambda $ up to its boundary, and
    there exists a universal constant $C $ such that 
    \begin{equation}\label{ee}
        \|\nabla e_j(t,\cdot)\|_{L^\infty(B_{1/2}\cap\Lambda )}\leq C\sqrt{K(0)}.
    \end{equation}
    Furthermore, $e_1(t,\cdot)$ is odd in the $x_1$-variable and $e_2(t,\cdot)$ is even in the $x_1$ variable within the lifespan of the solution, and they satisfy the following estimates 
    \[\exists C>0,~\forall x\in B_{1/2}\cap\Lambda ,~~|e_j(t,x)|\leq C\sqrt{K(0)}|x_j|,\quad j=1,2.\]
\end{lemma}
\begin{proof}
    A similar result was proved in \cite[Prop. 3.2-3.3]{HLY2024} for the case that the horizontal component lies in $\mathbb{T}$.
\end{proof}

\subsection{Estimates of velocity in terms of the vorticity}
Lemma \ref{le1} is an important pointwise velocity estimate generalized from \cite[Lemma 3.1]{Sverak2013}.
It is useful to control the velocity near the origin. Let $\Lambda_+=\RR_+\times[0,1] $.
\begin{lemma}\label{le1}
    Assume that $\omega(t,\cdot)\in (L^\infty \cap L^1)(\D_t)$ is the vorticity associated with a solution to the system \eqref{wave1}--\eqref{wave4}, 
    and that the initial data satisfy the assumptions of both Proposition \ref{pp2} and Lemma \ref{sym}.
    Then for any $x\in B_{1/2}\cap\Lambda_+ $, the following holds for all time during the lifespan of the solution: 
    \begin{equation}\label{le1-1}
        u_j(t,x)={(-1)}^j\left(\frac4\pi\int_{Q(|x|)}\frac{y_1 y_2}{{|y|}^4}\omega(t,y)\,dy+b_j(t,x)\right)x_j,
    \end{equation}
    where $Q(r)=\Lambda_+\setminus\overline{B_r(0)} $ and $b_1,b_2 $ satisfy
    \[|b_1(t,x)|\leq C_0\left(\|\omega\|_{L^\infty}\left(1+\log\left(1+\frac{x_2}{x_1}\right)\right)+\sqrt{K(0)}\right),\]
    \[|b_2(t,x)|\leq C_0\left(\|\omega\|_{L^\infty}\left(1+\log\left(1+\frac{x_1}{x_2}\right)\right)+\sqrt{K(0)}\right),\]
    for some universal constant $C_0 $.
\end{lemma}
\begin{proof}
    Recall that \eqref{error} reads
    \[u(t,x)=U(t,x)+e(t,x),\qquad x\in B_{1/2}\cap\Lambda_+,t\in[0,T),\]
    and by Lemma \ref{err}, we have 
    \[|e_j(t,x)|\leq C\sqrt{K(0)}|x_j|,\qquad j=1,2.\]
    Let $\tilde{\omega} $ be the odd-in-$x_2 $ extension of $\omega $ from $\Lambda=\RR\times[0,1] $ to $\RR\times[-1,1] $, namely
    \[\tilde{\omega}(t,x_1,x_2)=\begin{cases}
        \omega(t,x_1,x_2), & x_2\geq 0,\\
        -\omega(t,x_1,-x_2), & x_2<0.
    \end{cases}\]
    The elliptic equation of $\Psi $, namely \eqref{Psi}, admits a unique odd-odd solution given by the Biot-Savart law:
    \[\Psi(t,x)=\frac1{4\pi}\sum_{m\in\ZZ}\int_{\RR\times(-1,1)}\ln({(x_1-y_1)}^2+{(x_2-y_2+2m)}^2)\tilde{\omega}(t,y)\,dy,\]
    which leads to the following expression of $U$:
    \begin{align*}
        U(t,x)& =\nabla^{\bot}\Psi(t,x)\\
        & =\frac1{2\pi}\sum_{m\in\ZZ}\int_{\RR\times(0,1)}\left(\frac{r_m}{{|r_m|}^2}-\frac{s_m}{{|s_m|}^2}\right)\omega(t,y)\,dy\\
        & =\frac1{2\pi}\sum_{m\in\ZZ}\int_0^1\int_0^\infty \left(\frac{r_m}{{|r_m|}^2}-\frac{s_m}{{|s_m|}^2}-\frac{r_m'}{{|r_m'|}^2}+\frac{s_m'}{{|s_m'|}^2}\right)\omega(t,y_1,y_2)\,dy_1\,dy_2,
    \end{align*}
    where 
    \[r_m=(-y_2+x_2+2m,x_1-y_1),\quad s_m=(y_2+x_2+2m,x_1-y_1),\]
    \[r_m'=(-y_2+x_2+2m,x_1+y_1),\quad s_m'=(y_2+x_2+2m,x_1+y_1).\]
    
    For $m\neq 0 $, we estimate the integrand using $|r_m|,|r_m'|,|s_m|,|s_m'|\approx\max\{y_1,|m|\} $. Take $U_1 $ as an example, 
    \[\frac{-y_2+x_2+2m}{{|r_m|}^2}-\frac{-y_2+x_2+2m}{{|r_m'|}^2}=\frac{4 y_1 x_1(x_2-y_2+2m)}{{|r_m|}^2 {|r_m'|}^2}=x_1 O\left(\frac1{m^2}\right),\]
    and similarly,
    \[\frac{y_2+x_2+2m}{{|s_m|}^2}-\frac{y_2+x_2+2m}{{|s_m'|}^2}=\frac{4 y_1 x_1(x_2+y_2+2m)}{{|s_m|}^2 {|s_m'|}^2}=x_1 O\left(\frac1{m^2}\right).\]
    Analogous estimates hold for $U_2$. Therefore, 
    \[\left|\frac1{2\pi}\sum_{m\neq 0}\int_{(0,\infty)\times(0,1)}\left(\frac{r_m}{{|r_m|}^2}-\frac{s_m}{{|s_m|}^2}-\frac{r_m'}{{|r_m'|}^2}+\frac{s_m'}{{|s_m'|}^2}\right)\omega(t,y)\,dy\right|\leq C x_1\|\omega(t,\cdot)\|_{L^\infty}\]
    for some universal constant $C $.

    Thus, it remains to consider the term with $m=0 $ in the expression of $U $, namely the following integral:
    \[\frac2{\pi}\int_{(0,\infty)\times(0,1)}\left[\frac{x_1 y_1(x_2-y_2)}{{|x-y|}^2{|x-\tilde{y}|}^2}-\frac{x_1 y_1(x_2+y_2)}{{|x-\bar{y}|}^2{|x+y|}^2}\right]\omega(t,y)\,dy,\]
    where $\tilde{y}=(-y_1,y_2),\bar{y}=(y_1,-y_2) $.
    In such case, we can divide this integral into two parts $Q(2|x|) $ and $(0,\infty)\times(0,1)\setminus Q(2|x|) $.
    Since the detailed proof of $\eqref{le1-1}$ is the same as Zlato\v{s} \cite[Lemma 2.1]{Zlatos2025}, we omit the details.
\end{proof}

\section{Proof of the main theorem}
The proof is generalized from that of \cite[Theorem 1.1]{Zlatos2025}.
However, the situation here is more delicate due to the presence of the free surface.
Denote by $\D_t^+=\D_t\cap\{x_1> 0\}$ the right half of the fluid domain.
As in \cite{Zlatos2025}, for $\theta\in(0,e^{-1}) $, we define
\[\Omega_\theta=\{x\in \D_0^+:x_1\in (\theta,e^{-1})\mbox{ and }x_2<g(x_1)\},\quad\mbox{where }g(s)=se^{{|\ln s|}^{\frac12}}.\] 
The constant $\theta$ is understood as a spatial truncation parameter and the function $g$ is considered as an ``envelope function''.
Before proceeding to the proof of the main theorem, we shall establish the following lemma.
\begin{lemma}\label{le3-1}
    Let $\eps>0$ be the small constant in Theorem \ref{th1}. Under the assumptions of Proposition \ref{pp2}, there exists a constant $s_0\in(0,e^{-4})$ such that for any $\theta\in (0,s_0]$ and $\alpha\in(0,1/3]$, if $\omega\in (L^\infty\cap L^1)(\D_t) $ is odd in $x_1 $ and 
    $\varepsilon\chi_{\alpha\Omega_\theta}\leq \omega\chi_{\D_t^+}\leq \varepsilon\chi_{\D_t^+} $ with $\D_t^+=\D_t\cap\{x_1>0\} $, then the corresponding velocity field $u$ satisfies the following pointwise estimates on the boundary of $\alpha\Omega_\theta$:
    \begin{itemize}
 	\item [(1)] (Inner upper boundary, singularity-dominated) For any $s\in[\theta,s_0]$, the fluid velocity points strictly outward or is tangent to the boundary, that is,
   		\begin{equation}\label{le2-1}
            {(-1)}^j u_j(s\alpha,g(s)\alpha)\geq 0 \quad\quad \text{for } j=1,2.
        \end{equation}
	\item [(2)] (Outer upper boundary, error-dominated) For any $s \in [s_0, e^{-1}]$, the inward penetration of the velocity field possesses the following bound: There exists a constant $C'$ (independent of $\varepsilon,\alpha,\theta)$ such that
        \begin{equation}\label{le2-2}
            {(-1)}^j u_j(s\alpha,g(s)\alpha)\geq -C'\eps s\alpha  \quad\quad \text{for } j=1,2.
    	\end{equation}
    \item [(3)] (Left-vertical boundary) On the segment $\{x_1=\theta\alpha\}$, the horizontal velocity has the following upper bound
	    \begin{equation}\label{le2-3}
            \frac{1}{\theta \alpha} \sup_{s \in [0, g(\theta)\alpha]} u_1(\theta\alpha, s) \le \eps \Big( -h(\theta)|\ln \theta| + C_0 f(\theta) \Big) + \frac{e}{\alpha} \inf_{s \in [0, \alpha]} u_1(e^{-1}\alpha, s),
    	\end{equation}
        where $f(\theta) = 5+2\sqrt{C_1} + |\ln \theta|^{1/2}$ with $C_1 $ being the constant in Proposition \ref{pp2}, and $h(\theta)$ being a continuous function satisfying $\lim\limits_{\theta \to 0} h(\theta) = \frac{2}{\pi}$. 
    \end{itemize}
\end{lemma}

Lemma \ref{le3-1} gives the estimates of the velocity field in $\{(s\alpha,g(s)\alpha):s\in [\theta,e^{-1}]\} $.
It helps to control the flow of some part of the level set $\{\omega(t,\cdot)=\varepsilon\} $.

\begin{proof}
    From \cite[Section 3]{Zlatos2025}, we know that there exists some $h:(0,e^{-4})\to (0,\infty) $ with $\lim_{s\to 0}h(s) =2/\pi $ such that 
    \[\frac4{\pi}\int_{D_s}\frac{y_1 y_2}{{|y|}^4}\,dy=h(s)|\ln s|\quad s\in (0,e^{-4}),\]
    where 
    $D_s=\Omega_0\cap Q(\sqrt{s^2+{g(s)}^2})\cap B_{e^{-1}}(0) $.
    
    (1) Since $\alpha\Omega_\theta\subseteq  1/3\Omega_{s_0}\subseteq B_{1/2}\cap\Lambda_+$,
    we can use Lemma \ref{le1} (for $x=(s\alpha,g(s)\alpha) $). By $\lim\limits_{s\to 0}g(s)/s=\infty $, we obtain that: for sufficiently small $s_0>0 $ and any $s\in [\theta,s_0] $, there holds 
    \begin{align*}
        {(-1)}^j u_j(s\alpha,g(s)\alpha) & = \left(\frac4\pi\int_{Q(|x|)}\frac{y_1 y_2}{{|y|}^4}\omega(y)\,dy+b_j(x)\right)x_j\\
        & \geq \left(\varepsilon h(s)|\ln (s)|-C_0\varepsilon\left(1+\log\left(\frac{\alpha s+g(s)\alpha}{\alpha s}\right)+\sqrt{C_1}\right)\right)x_j\\
        & \geq \left(\varepsilon h(s)|\ln (s)|- C_0\varepsilon\left(2+{|\ln s|}^{\frac12}+\sqrt{C_1}\right)\right)x_j \geq 0.
    \end{align*} 
    
    (2) For $s\in [s_0,e^{-1}] $, we have 
    \begin{align*}
        {(-1)}^j u_j(s\alpha,g(s)\alpha) & = \left(\frac4\pi\int_{Q(|x|)}\frac{y_1 y_2}{{|y|}^4}\omega(y)\,dy+b_j(x)\right)x_j\\
        & \geq \left(-C_0\varepsilon\left(1+\log\left(\frac{\alpha s+g(s)\alpha}{\alpha s}\right)+\sqrt{C_1}\right)\right)x_j\\
        & \geq -C'\eps\frac{s}{g(s)}x_j\geq -C'\varepsilon\alpha s,
    \end{align*}
    where 
    \[C'=\max_{s\in[s_0,e^{-1}]}C_0 \frac{g(s)}{s}\left(1+\sqrt{C_1}+\ln\frac{s+g(s)}{s}\right).\]
    
    (3) Let $x=(\theta\alpha,g(\theta)\alpha) $. Then, using Lemma \ref{le1}, we obtain
    \begin{align*}
        \frac1{(\theta\alpha)}\sup_{s \in [0, g(\theta)\alpha]}u_1(\theta\alpha,s)=&\sup_{s \in [0, g(\theta)\alpha]}-\left(\frac4{\pi}\int_{Q(|(\theta\alpha,s)|)}\frac{y_1 y_2}{{|y|}^4}\omega(y)\,dy+b_1(\theta\alpha,s)\right)\\
        \leq & -\frac4{\pi}\int_{Q(|x|)}\frac{y_1 y_2}{{|y|}^4}\omega(y)\,dy+C_0\varepsilon\left(1+\log\left(\frac{\theta\alpha+g(\theta)\alpha}{\theta\alpha}\right)+\sqrt{C_1}\right)\\
        \leq & -\frac4{\pi}\int_{Q(e^{-1}\alpha)}\frac{y_1 y_2}{{|y|}^4}\omega(y)\,dy-\varepsilon h(\theta)|\ln\theta|+C_0\varepsilon\left(2+{|\ln \theta|}^{\frac12}+\sqrt{C_1}\right)
    \end{align*}
    and
    \begin{align*}
        \frac{e}{\alpha} \inf_{s \in [0, \alpha]} u_1(e^{-1}\alpha, s)=& \inf_{s \in [0, \alpha]}-\left(\frac4{\pi}\int_{Q(|(e^{-1}\alpha,s)|)}\frac{y_1 y_2}{{|y|}^4}\omega(y)\,dy+b_1(e^{-1}\alpha,s)\right)\\
        \geq & -\frac4{\pi}\int_{Q(e^{-1}\alpha)}\frac{y_1 y_2}{{|y|}^4}\omega(y)\,dy-C_0\varepsilon\left(1+\log\left(1+\frac{\alpha}{e^{-1}\alpha}\right)+\sqrt{C_1}\right)\\
        \geq & -\frac4{\pi}\int_{Q(e^{-1}\alpha)}\frac{y_1 y_2}{{|y|}^4}\omega(y)\,dy-C_0\varepsilon\left(3+\sqrt{C_1}\right).
    \end{align*}
    Combining the above estimates, we conclude that \eqref{le2-3} holds.
\end{proof}

Now we are ready to prove Theorem \ref{th1}.
\begin{proof}[Proof of  Theorem \ref{th1}]
Let $\theta:[0,\infty)\to (0,s_0] $ be a decreasing function (to be determined) and let $\omega $ be an odd (in $x_1 $) solution such that $\omega(0,\cdot)\in C^\infty(\D_0) $ and
$\varepsilon\chi_{\frac13\Omega_{\theta(0)}}\leq \omega(0,\cdot)\chi_{\D_0^+}\leq\varepsilon\chi_{\D_0^+} $.
Also we assume $\|\omega(0,\cdot)\|_{L^\infty}=\varepsilon $, which then implies $\|\omega(t,\cdot)\|_{L^\infty}=\varepsilon $.
So we can find an initial velocity $u_0 $ and a constant $C_1 $ such that $K(0)=C_1\varepsilon^2 $.

\textbf{Step 1: Selecting the shrinking scale factor $\alpha(t)$.} 
We define $\alpha:[0,\infty)\to (0,1/3] $ to be the function satisfying $\alpha(0)=1/3 $ and 
\begin{equation}\label{d}
    \alpha'(t)=-\eps\left(2 C' e^{-1} \rho_0^{-1}+C_0\big(3+\sqrt{C_1}\big)\right)\alpha(t)+e\inf_{s\in[0,\alpha(t)]}u_1(t,e^{-1}\alpha(t),s),\quad t\geq 0,
\end{equation}
where 
\[\rho_0= \frac{s_0}{2{|\ln s_0|}^{1/2}}\leq \min_{s\in [s_0,e^{-1}]}\left(\frac{g(s)}{g'(s)}-s\right)\] is strictly positive due to the strict convexity of $g$. The function $\alpha $ represents the shrinking rate of the region in $\{\omega(t,\cdot)=\varepsilon\} $.
Furthermore, combining this with Lemma \ref{le1}, we have 
\[u_1(t,e^{-1}\alpha,s)\leq C_0\varepsilon\left(1+\log\left(1+\frac{\alpha}{e^{-1}\alpha}\right)+\sqrt{C_1}\right)e^{-1}\alpha\leq C_0\varepsilon\left(3+\sqrt{C_1}\right)e^{-1}\alpha\]
for any $s\in [0,\alpha] $. This estimate and \eqref{d} imply that 
\begin{equation}\label{d2}
    \alpha'(t)\leq -2C'\eps e^{-1}\rho_0^{-1}\alpha(t).
\end{equation}

\textbf{Step 2: Verifying the Non-penetration.}  We now want to show that, if we appropriately select the function $\theta(t) $, then the boundary of the level set $\{\omega(t,\cdot)=\varepsilon\}$ never enters $\alpha(t)\Omega_{\theta(t)}$. 
The function $\theta(t) $ describes how fast this region propels to the origin. To ensure this, the boundary of $\alpha(t)\Omega_{\theta(t)}$ must shrink inward faster than the ``shrinking speed'' of the fluid, and its left boundary must push inward more slowly than the ``inward fluid velocity''. We verify this on all three boundary segments.

\textbf{Part 1: Right vertical boundary $\{x_1=e^{-1}\alpha(t)\}$.} Indeed, by $\eqref{d}$, we have 
\[u_1(t,e^{-1}\alpha(t),s)>e^{-1}\alpha'(t)\quad\text{ on }\{(e^{-1}\alpha(t),s):s\in [0,\alpha(t)]\}.\]

\textbf{Part 2: Upper curved boundary.} This part can be parametrized as $\{(s\alpha(t),g(s)\alpha(t)):s\in[\theta(t),e^{-1}]\} $. 
The outward normal is proportional to $n = (-1, {g'(s)}^{-1})$, and therefore the normal projection of the boundary's moving velocity $V_0 $ is 
\[V_0\cdot n = \alpha'(t) \left( \frac{g(s)}{g'(s)} - s \right).\]
The normal component of the fluid velocity is $u \cdot n = -u_1 + u_2 / g'(s)$. We analyze this in two sub-regions.
\begin{itemize}
\item For the inner region $s \in [\theta(t), s_0]$: Lemma \ref{le3-1}(1) gives $-u_1 \geq 0$ and $u_2 \geq 0$, thus $u \cdot n \geq 0$. Since $\alpha'(t) < 0$ and $(g(s)/g'(s) - s) > 0$, we have $V_0 \cdot n < 0 \leq u \cdot n$. The condition is strictly satisfied.
\item For the outer region $s \in [s_0, e^{-1}]$: Lemma \ref{le3-1}(2) gives  the bounds $-u_1 \geq -C'\eps s\alpha$ and $u_2 \geq -C'\eps s\alpha$. Thus, 
\[u \cdot n \ge -C'\varepsilon s\alpha - \frac{C'\varepsilon s\alpha}{g'(s)}\ge -2 C's\varepsilon\alpha(t) \ge -2 C' e^{-1}\varepsilon\alpha(t).\] 
To ensure $u \cdot n \geq V_0 \cdot n$, it suffices to show $-2C'e^{-1}\varepsilon\alpha(t) \ge \alpha'(t) \rho_0$, which is equivalent to $\alpha'(t)\leq -2C'\eps e^{-1}\rho_0^{-1}\alpha(t)$. 
Due to \eqref{d2}, this inequality definitely holds.
\end{itemize}

\textbf{Part 3: Left vertical boundary.}  
This is the core region driving the double-exponential growth. The fluid experiences strong convective transport toward the left ($u_1 < 0$). 
We shall prove that the left boundary's leftward speed must be slower than the fluid's leftward speed.  The boundary's normal speed is $-(\theta(t)\alpha(t))'$. The fluid normal velocity is $-u_1$. 
We shall require $-(\theta\alpha)' \le -u_1$, which can also be written as 
\[\dfrac{\theta'}{\theta} + \dfrac{\alpha'}{\alpha} \geq \dfrac{u_1}{\theta\alpha}.\]
In order to establish this inequality for all points on the left boundary $\{x_1 = \theta(t)\alpha(t)\}$, the left side must be $\geq$ the supremum of the right side. Using Lemma \ref{le3-1}(3), we get 
\[\frac{\theta'}{\theta} + \frac{\alpha'}{\alpha} \geq \eps \left( -h(\theta)|\ln \theta| + C_0 f(\theta) \right) + \frac{e}{\alpha} \inf_{s\in[0,\alpha]} u_1(e^{-1}\alpha, s).\] 
Invoking the expression of $\alpha'/\alpha$, we find a cancellation, which yields that  $\theta(t)$ must satisfy
\begin{align}
\label{imp} \frac{\theta'(t)}{\theta(t)} \geq \eps \left[ -h(\theta(t))|\ln \theta(t)| + C_0 f(\theta(t)) + \left(2 C' e^{-1} \rho_0^{-1}+C_0\big(3+\sqrt{C_1}\big)\right) \right].
\end{align}

In other words, we prove that the boundary of $\{\omega(t,\cdot)=\varepsilon\} $ never enters $\alpha(t)\Omega_{\theta(t)} $, provided that $\eqref{imp} $ holds true. 
This will then lead to $\omega(t,\theta(t)\alpha(t),s)=\varepsilon $ for all $t\geq 0 $ and $s\in(0,g(\theta(t))\alpha(t)]$.  Now, since we have $\omega(t,0,s)=0 $, we then get the following inequality
\[\|\nabla\omega(t,\cdot)\|_{L^\infty(\D_t)}\geq \frac{\varepsilon}{\theta(t)\alpha(t)}\geq \frac{\varepsilon}{\theta(t)}.\]

\textbf{Step 3: Solving the differential inequality for $\theta(t)$.}
In fact, we can choose the function $\theta $ such that it satisfies  
\[\lim_{t\to\infty}\frac{\ln(-\ln\theta(t))}{t}=c_2\varepsilon\]
for some constant $c_2>0 $.
For example, we choose $\theta(t)=\exp(-\exp(\eta(t))) $, where $\eta(t)=\frac1{\pi}\varepsilon t+\ln c_1$ and $c_1 $ is a sufficiently large constant. Then a straightforward calculation shows that
\[\eta'(t) \leq \eps h(\theta(t)) - \eps e^{-\eta(t)} \left[ C_0 f(\theta(t)) + \left(2 C' e^{-1} \rho_0^{-1}+C_0\big(3+\sqrt{C_1}\big)\right) \right],\]
where $f(\theta)=5+2\sqrt{C_1}+|\ln \theta|^{1/2} $.
Thus, we obtain the double exponential growth for the vorticity gradient: for all $t\geq 0 $ within the lifespan of the solution, we have 
\[\|\nabla\omega(t,\cdot)\|_{L^\infty(\D_t)}\geq\frac{\varepsilon}{\theta(t)}\geq\varepsilon\exp\left(c_1\exp\left(\frac{\varepsilon t}{\pi}\right)\right).\]
\end{proof}


\begin{thebibliography}{99}
\bibitem{ABZ2011}
Thomas Alazard, Nicolas Burq, Claude Zuily.
\newblock On the water-wave equations with surface tension.
\newblock {\em Duke Math. J.}, 158(3):413--499, 2011.

\bibitem{ABZ2014}
Thomas Alazard, Nicolas Burq, Claude Zuily.
\newblock On the Cauchy problem for gravity water waves.
\newblock {\em Invent. Math.}, 198(1):71--163, 2014.

\bibitem{AD2015}
Thomas Alazard, Jean-Marc Delort.
\newblock Global solutions and asymptotic behavior for two-dimensional gravity water waves.
\newblock {\em Ann. Sci. {\'E}c. Norm. Sup{\'e}r.(4)}, 48(5):1149--1238, 2015.


\bibitem{Steiner}
Ennio De Giorgi.
\newblock Sulla propriet\`a isoperimetrica dell'ipersfera, nella classe degli insiemi aventi frontiera orientata di misura finita.
\newblock {\em Atti Accad. Naz. Lincei Mem. Cl. Sci. Fis. Mat. Nat. Sez. I (8)}, 5:33--44, 1958.

\bibitem{CS2007}
Daniel Coutand, Steve Shkoller.
\newblock Well-posedness of the free-surface incompressible Euler equations with or without surface tension.
\newblock {\em J. Amer. Math. Soc.}, 20(3):829--930, 2007.
  
\bibitem{DIPP2017}
Yu Deng, Alexandru D. Ionescu, Beno\^{i}t Pausader, Fabio Pusateri.
\newblock Global solutions of the gravity-capillary water-wave system in three dimensions.
\newblock {\em Acta Math.}, 219(2):213--402, 2017.

\bibitem{Denisov2009}
Sergey Denisov, 
\newblock Infinite superlinear growth of the gradient for the two-dimensional {Euler} equation.
\newblock {\em Discrete Contin. Dyn. Syst.}, 23 (2009), 755--764.


\bibitem{GMS2012}
Pierre Germain, Nader Masmoudi, Jalal Shatah.
\newblock Global solutions for the gravity water waves equation in dimension 3.
\newblock {\em Ann. Math.}, 175(2), 691--754, 2012.

\bibitem{GMS2015}
Pierre Germain, Nader Masmoudi, Jalal Shatah.
\newblock Global existence for capillary water waves.
\newblock {\em Commun. Pure. Appl. Math.}, 68(4):625--687, 2015.

\bibitem{HLY2024}
Zhongtian Hu, Chenyun Luo, and Yao Yao.
\newblock Small scale creation for 2D free boundary Euler equations with
  surface tension.
\newblock {\em Ann. PDE}, 10(2), 2024.

\bibitem{IT2016}
Mihaela Ifrim, Daniel Tataru.
\newblock Two-dimensional water waves in holomorphic coordinates II: global solutions.
\newblock {\em  Bull. Soc. Math. France}, 144(2):366--394, 2016.

\bibitem{IT2017}
Mihaela Ifrim, Daniel Tataru.
\newblock The lifespan of small data solutions in two-dimensional capillary water waves.
\newblock {\em Arch. Rational Mech. Anal.}, 225(3):1279--1346, 2017.

\bibitem{IP2015}
Alexandru D. Ionescu, Fabio Pusateri.
\newblock Global solutions for the gravity water waves system in 2D.
\newblock {\em Invent. Math.}, 199(3):653--804, 2015.

\bibitem{IP2018}
Alexandru D. Ionescu, Fabio Pusateri.
\newblock Global regularity for 2D water waves with surface tension.
\newblock {\em Mem. Amer. Math. Soc.}, 256(1227), 2018.

\bibitem{Sverak2013}
Alexander Kiselev, Vladimir {\v{S}}ver{\'a}k.
\newblock Small scale creation for solutions of the incompressible two-dimensional Euler equation.
\newblock {\em Ann. Math.}, 180(3), 1205--1220, 2014.

\bibitem{JYZ2025}
In-Jee Jeong, Yao Yao, Tao Zhou.
\newblock Superlinear gradient growth for 2D Euler equation without boundary.
\newblock arXiv:2507.15739, preprint, 2025.

\bibitem{Lannes2005}
David Lannes.
\newblock Well-posedness of the water-waves equations.
\newblock {\em J. Amer. Math. Soc.}, 18(3):605--654, 2005.

\bibitem{Lindblad2005}
Hans Lindblad.
\newblock Well-posedness for the motion of an incompressible liquid with free surface boundary.
\newblock {\em Ann. Math.}, 162(1), 109--194, 2005.

\bibitem{MZ2009}
Mei Ming, Zhifei Zhang.
\newblock Well-posedness of the water-wave problem with surface tension.
\newblock {\em J. Math. Pures Appl.}, 92(5), 429--455, 2009.

\bibitem{SZ1}
Jalal Shatah, Chongchun Zeng.
\newblock Geometry and a priori estimates for free boundary problems of the {E}uler's equation.
\newblock {\em Commun. Pure. Appl. Math.}, 61(5), 698--744, 2008.

\bibitem{SZ2}
Jalal Shatah, Chongchun Zeng.
\newblock A priori estimates for fluid interface problems.
\newblock {\em Commun. Pure. Appl. Math.}, 61(6), 848--876, 2008.

\bibitem{SZ3}
Jalal Shatah, Chongchun Zeng.
\newblock Local well-posedness for fluid interface problems.
\newblock {\em Arch. Rational Mech. Anal.}, 199(2), 653--705, 2011.
  
\bibitem{Su2020}
Qingtang Su.
\newblock Long time behavior of 2D water waves with point vortices.
\newblock {\em Commun. Math. Phys.}, 380(3):1173--1266, 2020.

\bibitem{WangXC2018}
Xuecheng Wang.
\newblock Global infinite energy solutions for the 2D gravity water waves system.
\newblock {\em Commun. Pure. Appl. Math.}, 71(1):90--162, 2018.

\bibitem{WangXC2020}
Xuecheng Wang.
\newblock Global regularity for the 3D finite depth capillary water waves.
\newblock {\em Ann. Sci. Éc. Norm. Supér.}, 53(4):847--893, 2020.

\bibitem{Wu1997}
Sijue Wu.
\newblock Well-posedness in {S}obolev spaces of the full water wave problem in 2-{D}.
\newblock {\em Invent. Math.}, 130(1):39--72, 1997.

\bibitem{Wu1999}
Sijue Wu.
\newblock Well-posedness in {S}obolev spaces of the full water wave problem in 3-{D}.
\newblock {\em J. Amer. Math. Soc.}, 12(2):445--495, 1999.

\bibitem{Wu2009}
Sijue Wu.
\newblock Almost global wellposedness of the 2-D full water wave problem.
\newblock {\em Invent. Math.}, 177(1):45--135, 2009.

\bibitem{Wu2011}
Sijue Wu.
\newblock Global wellposedness of the 3-D full water wave problem.
\newblock {\em Invent. Math.}, 184(1):125--220, 2011.

\bibitem{Xu2016}
Xiaoqian Xu.
\newblock Fast growth of the vorticity gradient in symmetric smooth domains for 2{D} incompressible ideal flow.
\newblock {\em J. Math. Anal. Appl.}, 439(2), 594--607, 2016.


\bibitem{ZZ2008}
Ping Zhang, Zhifei Zhang.
\newblock On the free boundary problem of three-dimensional incompressible {E}uler equations.
\newblock {\em Commun. Pure. Appl. Math.}, 61(7):877--940, 2008.

\bibitem{Zlatos2015}
Andrej Zlato\v{s}.
\newblock Exponential growth of the vorticity gradient for the Euler equation on the torus.
\newblock {\em Adv. Math.}, 268:396--403, 2015.

\bibitem{Zlatos2025}
Andrej Zlato\v{s}.
\newblock Maximal double-exponential growth for the Euler equation on the half-plane.
\newblock {\em Invent. Math.}, 243(1):117--126, 2026.

\bibitem{ZhengF2019}
Fan Zheng.
\newblock Long-term regularity of 3D gravity water waves.
\newblock {\em Commun. Pure Appl. Math.}, 75(5):1074--1180, 2022.

\end{thebibliography}
\end{document}